# A Tight Configuration-Component Based Hybrid Model for Combined-Cycle Units in MISO Day-Ahead Market

Chenxi Dai, *Student Member, IEEE*, Yonghong Chen, *Senior Member, IEEE*, Fengyu Wang, *Member, IEEE*, Jie Wan, *Member, IEEE*, and Lei Wu, *Senior Member, IEEE*

***Abstract*—This paper proposes a hybrid CCGT model for day-ahead market clearing, in order to enhance operation flexibility of CCGTs in practice. The proposed hybrid model, by taking benefits of combined offers on configurations and individual PTs from MPs, can reflect physical operation features of CCGTs more accurately. A comprehensive review on CCGT models in literature is conducted, and their advantages and shortcomings are analyzed. The mapping relationship between the two most investigated models, i.e., CFBM and CPBM, is then revealed and constructed for deriving the proposed hybrid model. Tightening formulations are further discussed for achieving better computational performance. The proposed hybrid model is tested and compared with CFBM via an IEEE 118-bus system and the MISO system. Results show great benefits in maintaining operation flexibility and enhancing social welfare.**

## NOMENCLATURE

Major symbols and notations used throughout the paper are listed below, while others are defined following their first appearances as needed.

**Acronyms:**

| | |
|---|---|
| CAISO | California Independent System Operator |
| CCGT | Combined-cycle gas turbine |
| CFBM | Configuration-based model |
| CPBM | Component-based model |
| CT/ST/DB | Combustion turbine/steam turbine/duct burner |
| ERCOT | Electric Reliability Council of Texas |
| HRSG | Heat recovery steam generator |
| ISO-NE | Independent System Operator New England |
| LHS | Left-hand-side |
| MISO | Midcontinent Independent System Operator |
| MP | Market participant |
| PT | Physical turbine, including CT and ST |
| SCUC | Security-constrained unit commitment |
| SPP | Southwest Power Pool |

**Indices:**

| | |
|---|---|
| $c, s$ | Indices of CTs and STs per CCGT, i.e., $c, s \in X$ |
| $j$ | Index of CCGTs, i.e., $j \in J$ |
| $t, t'$ | Indices of time intervals (hours), i.e., $t, t' \in T$ |
| $w$ | Type of PT startup warmth status, i.e., $w \in W$ |
| $x$ | Index of PTs per CCGT, $x \in X$ |
| $y, y'$ | Indices of configurations per CCGT, i.e., $y, y' \in Y$ |

**Sets:**

| | |
|---|---|
| $F_j^{yy'}$ | Set of feasible transitions of CCGT $j$, where $y/y'$ is the "from"/"to" configuration and $y \neq y'$ |
| $J/\ T$ | Set of CCGTs/ time intervals |
| $ON^x$ | Set of valid configurations with PT $x$ on |
| $UT^x$,$DT^x$ | Set of upward/downward transitions which include the startup/shutdown of PT $x$ |
| $W$ | Set of startup type warmth status, i.e., $W$={1=hot; 2=intermediate; 3=cold} |
| $X_j/\ Y_j$ | Set of PTs/valid configurations of CCGT $j$ |

**Variables:**

| | |
|---|---|
| $i_{j,t}^x$ | Intermediate variable for commitment status of PT $x$ for CCGT $j$ at $t$ |
| $p_{j,t}^y/u_{j,t}^y$ | Power output/commitment status of CCGT $j$ at configuration $y$ and $t$ |
| $su_{j,t}^x/sd_{j,t}^x$ | Intermediate variable for startup/ shutdown indicators of PT $x$ for CCGT $j$ at $t$ |
| $v_{j,t}^{yy'}$ | Transition from $y$ to $y'$ ($y \neq y'$) for CCGT $j$ at $t$ |
| $\delta_{j,t,w}^x$ | Indicator of startup type $w$ of PT $x$ for CCGT $j$ at time $t$ |

**Constants:**

| | |
|---|---|
| $C_j^{NL,y}$ | No-load cost of CCGT $j$ at configuration $y$ |
| $C_{j,w}^{SU,x}$ | Startup cost of PT $x$ for CCGT $j$ with type $w$ |
| $C_j^{P,y}(\cdot)$ | Incremental energy cost function for CCGT $j$ at configuration $y$ |
| $i_j^{0,x}$ | Initial commitment status of PT $x$ for CCGT $j$ |
| $\underline{P}_{j,t}^y/\overline{P}_{j,t}^y$ | Generation resource minimum/maximum output limit for CCGT $j$ at configuration $y$ and time $t$ |
| $RU_j^y/RD_j^y$ | Ramp up/down rate for configuration $y$ of CCGT $j$ at non-shutdown intervals |
| $RU0_j^y$ | Ramp up rate for configuration $y$ of CCGT $j$ at startup intervals |
| $T_w^{SU,x}$ | Minimum time that PT $x$ must be down for startup type $w$ |
| $TD_j^x/TU_j^x$ | Minimum down/up time of PT $x$ for CCGT $j$ |
| $TD_j^{0,x}/TU_j^{0,x}$ | Number of time intervals that PT $x$ of CCGT $j$ has been off /on before the scheduling horizon |

## I. INTRODUCTION

### A. Motivation

As the electricity grid accounts for 35% of U.S. total $CO_2$ emissions in 2016 as the largest contributor, more stringent environmental regulations on power plants become the most critical concern. Indeed, according to the Clean Power Plan 2015, the U.S. power sector targets to reduce 32% from its 2005 $CO_2$ emission level by 2030 [1]-[2]. To achieve

C. Dai and L. Wu are with Electrical and Computer Engineering Department, Clarkson University, Potsdam, NY 13699 USA. (e-mail: daic@clarkson.edu, lwu@clarkson.edu). Y. Chen and F. Wang are with Midcontinent Independent System Operator, Inc. (MISO), Carmel, IN, 46032 USA. J. Wan is with GE Grid Software Solution, Redmond, WA, 98052, USA. Disclaimer: The views expressed in this paper are solely those of the authors and do not necessarily represent those of MISO or GE Grid Software Solution.

this goal, one effective building block is to increase power generation from lower-emission CCGTs.

Compared to traditional coal-fired units, CCGTs enjoy many advantages including relatively shorter installation time, lower levelized cost of electricity, faster response time, less $CO_2$ emissions, as well as higher operation efficiency and reliability [3]-[5]. Indeed, different from old-fashioned coal-fired units that are incapable of providing sufficient flexibility, faster ramping capability of CCGTs makes them ideal to support deeper penetration of highly variable and uncertain renewable resources such as wind and solar. In addition, because of the new horizontal drilling and hydraulic fracturing techniques, huge shale gas reserves have been discovered in the U.S. and natural gas becomes commercially sustainable with a price comparable to that of coal. Therefore, CCGTs have made up the majority of generation capacities developed by private investors in the last decade, and an even higher proliferation of CCGTs is expected in the close future [6]. For instance, in ISO-NE, roughly 3,000 MW of new CCGTs is coming forward in recent auctions and will be available in 2017-2019 [7].

However, due to complicated operation characteristics, CCGTs bring significant challenges to the modeling and market clearing processes. In fact, for a CCGT, exhausted gas from CTs can be used by HRSG to produce steam and drive STs, while switching on/off CTs and STs needs to follow a required sequence. Fig. 1 shows allowed configurations and feasible transitions for a CCGT with 2 CTs and 1 ST (denoted as "2 × 1"). In addition, individual configurations corresponding to different CT/ST combinations may present different physical features, e.g., ramp rates, operation cost, etc.

Managing one of the largest electricity markets in the world, MISO encountered computational challenges in solving day-ahead SCUC in early 2014 due to significant expansion of its market footprint. Great efforts have been devoted since then to improve computational performance [8]-[10]. Recently, with the upward CCGT installation trend, MISO received requests from MPs to improve modelling capability. Better CCGT modelling in the market clearing engine can reduce production costs and incentivize MPs to offer costs closer to reality. In response to the request from MPs and the needs to improve market efficiency, a new CCGT model is proposed in this paper.

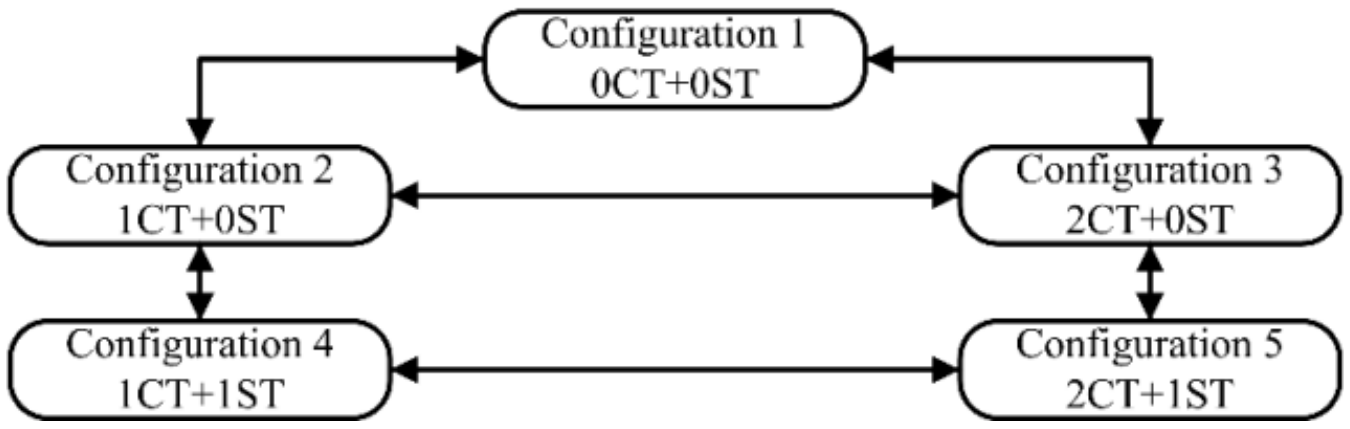

Fig.1. Feasible Transitions of a CCGT with 2 CTs and 1 ST

## B. *Literature Review*

As summarized in Table I, five CCGT models are currently used in academic research and industry practice.

TABLE I SUMMARY OF CCGT MODELS IN LITERATURE

| Model | 1 | 2 | 3 | 4 | 5 |
|---|---|---|---|---|---|
| Literature | N/A | [11] | [5], [14] | [4], [5], [10], [12], [14]-[23] | [24] |
| Industry Application | MISO, PJM, NYISO, XM | ISO-NE | N/A | ERCOT, SPP, CAISO, MISO prototype | N/A |

1) Aggregate Model

Aggregate model simplifies each CCGT as a regular unit by allowing the CCGT owner to provide an aggregate offer for the entire plant. This model sacrifices the accuracy because the aggregate offer cannot rigorously represent operating costs of CCGTs [11] or physical operating characteristics of individual PTs. As a result, SCUC engine may dispatch CCGTs at output levels that the PTs cannot achieve [12].

2) Pseudo Unit Model

Pseudo Unit Model represents a CCGT via multiple pseudo units, each of which comprises a single CT and its associated portion of the ST capacity. To ensure consistent commitment statuses of a same ST in multiple pseudo units, all pseudo units are required to have the same characteristics, including startup cost, no-load cost, incremental energy cost, minimum/maximum capacities, and minimum on/off times. Consequently, the pseudo unit model fails to precisely represent these physical characteristics, which are usually different under distinct operation configurations.

3) Component-Based Model (CPBM)

CPBM simulates individual components (i.e., PTs) of CCGTs separately, which can preserve generation characteristics of PTs and accurately reflect certain important features such as minimum up/down time and startup/shutdown cost. In addition, each PT has a market-clearing dispatch value that can be used for power flow assessment and contingency analysis [14]. However, this model has not been applied in practice due to the following issues: (i) with the unfriendly offer format for MPs [5], feasible transitions are hard to describe; and (ii) due to the strong dependency between CTs and STs, certain characteristics of individual PTs provided by MPs, such as energy costs and ramp rates, are not simply additive from PT offers.

4) Configuration-Based Model (CFBM)

In CFBM, different combinations of CTs/STs are aggregated into mutually exclusive configurations, as shown in Fig. 1. It is the most investigated model, which can better represent operational dependency between CTs and STs of a CCGT, especially incremental energy costs, ramp rates, and feasible transitions. Additionally, in order to mitigate the increased computational complexity and challenges to the market clearing process, practically some ISOs/RTOs limit the number of registered configurations per CCGT to reduce computational burden [16]. MISO has been evaluating CFBM for several years [9], [15], and benefit study shows that it could potentially bring over $16 million annual cost saving over the aggregate model [10]. With the dedicated computational research, MISO has greatly improved its market clearing software performance with the CFBM modeling [10] and starts to consider the implementation [9].

However, treating each configuration as a regular unit while not distinguishing individual PTs underneath each configuration, CFBM may undermine the accuracy of certain PT operation characteristics such as min up/down times and time-dependent transition costs. In practice, this could result in inefficiency/infeasibility. To better capture min up time limits of individual PTs, SPP introduced "group run time constraints" [17], in which a group represents the set of configurations that share a common operational constraint in a CCGT. In addition, the transition cost is usually derived from PT startup costs [10], while approximating warmth startup type of each configuration via CCGT "offline" status [18].

5) Edge-Based Model

Noticing the issues from previous models, edge-based model was proposed in [24] to capture physical features of individual PTs such as min up/down constraints. However, the edge-based model uses an entirely different logic to represent CCGT operations, which deviates prevailing modeling strategies in practice. This could increase additional logical complexity to formulate the problem and bring unnecessary efforts for practical implementation.

### *C. Contribution and Paper Organization*

In order to represent physical operation features of CCGTs in day-ahead market more accurately while alleviating computational burden, we extend the model in our existing work [25] to propose a new hybrid CCGT model in this paper. The main contributions of the paper include following aspects.

1) Driven by the needs from MPs, a hybrid CCGT model is proposed by taking advantages of offers from configure-tions and individual PTs, which would reflect CCGT features more accurately, provide more operation flexibility, and avoid infeasibility of PT operations.
2) Based on the proposed hybrid CCGT model, tightening formulations and techniques are further conducted to relieve computational burden. The proposed hybrid model is also compared with CFBM in terms of the size of formulation and potential computational impacts.
3) Test results on an IEEE 118-bus system and the MISO system show the proposed hybrid model can achieve more cost saving with reasonable computational performance compared with tightened CFBM formulation.

The remaining paper is organized as follows. Section II presents the proposed hybrid CCGT formulation. Section III discusses efforts for achieving tighter formulation and better computational performance. Case studies are shown in Section IV, and conclusion is drawn in Section V.

## II. Mathematical Formulation of the Hybrid Model

The mappings between configurations and PTs of a CCGT are first revealed, for both complete and aggregate configuration models. Various constraints are then formulated via configuration or individual PT offers accordingly, which have been confirmed by MPs that hybrid offers could better describe CCGT operation features. For the sake of concise representation, subscript $j$ for CCGTs is omitted in this section, and full SCUC formulation is not presented which can be referred in authors' work [8] [26].

### *A. Mapping Between Configurations and Individual PTs*

1) Complete Configuration Model

Using the "2×1" example as in Fig. 1, we illustrate the mapping relationship between configurations and individual PTs of a CCGT in the proposed hybrid model. When 2 CTs are considered non-identical, i.e., with different incremental energy costs, ramp rates, minimum up/down times or capacities, 7 configurations are required for representing all valid operating modes. Those are all off, $CT_1$, $CT_2$, $CT_1$+$CT_2$, $CT_1$+ST, $CT_2$+ST, and $CT_1$+$CT_2$+ST, and their commitment statuses, denoted as $u_t^1$ to $u_t^7$, are mutually exclusive at each time interval. In addition, 10 upward transitions between valid configurations in successive hours, with startup of certain PTs, are feasible, denoted as $v_t^{12}, v_t^{13}, v_t^{14}, v_t^{24}, v_t^{25}, v_t^{34}, v_t^{36}, v_t^{47}$, $v_t^{57}$, and $v_t^{67}$. Similarly, 10 downward transitions with shutdown of certain PTs are feasible, and reverse superscripts are used to represent downward transitions. Finally, the mapping relationship between commitment statuses of PTs $i_t^x$ and configurations $u_t^y$ is expressed as in (1), where $\boldsymbol{M}$ denotes the mapping relationship between PTs and configurations.

$$\boldsymbol{I} = \begin{bmatrix} i_t^{c1} \\ i_t^{c2} \\ i_t^{s1} \end{bmatrix} = \begin{bmatrix} 0 & 1 & 0 & 1 & 1 & 0 & 1 \\ 0 & 0 & 1 & 1 & 0 & 1 & 1 \\ 0 & 0 & 0 & 0 & 1 & 1 & 1 \end{bmatrix} \begin{bmatrix} u_t^1 \\ \vdots \\ u_t^7 \end{bmatrix} = \boldsymbol{M} \cdot \boldsymbol{U} \tag{1}$$

Similarly, (2)/(3) describes the relationship between PTs' startup/shutdown variables and feasible upward/downward transitions, where $\boldsymbol{M^{up}}/\boldsymbol{M^{dn}}$ represents mapping between start up/shutdown of PTs and feasible transitions of configurations.

$$\boldsymbol{SU} = \begin{bmatrix} su_t^{c1} \\ su_t^{c2} \\ su_t^{s1} \end{bmatrix} = \begin{bmatrix} 1 & 0 & 1 & 0 & 0 & 1 & 0 & 0 & 0 & 1 \\ 0 & 1 & 1 & 1 & 0 & 0 & 0 & 0 & 1 & 0 \\ 0 & 0 & 0 & 0 & 1 & 0 & 1 & 1 & 0 & 0 \end{bmatrix} \begin{bmatrix} v_t^{12} \\ \vdots \\ v_t^{67} \end{bmatrix} = \boldsymbol{M^{up}} \cdot \boldsymbol{V^{up}} \tag{2}$$

$$\boldsymbol{SD} = \begin{bmatrix} sd_t^{c1} \\ sd_t^{c2} \\ sd_t^{s1} \end{bmatrix} = \begin{bmatrix} 1 & 0 & 1 & 0 & 0 & 1 & 0 & 0 & 0 & 1 \\ 0 & 1 & 1 & 1 & 0 & 0 & 0 & 0 & 1 & 0 \\ 0 & 0 & 0 & 0 & 1 & 0 & 1 & 1 & 0 & 0 \end{bmatrix} \begin{bmatrix} v_t^{21} \\ \vdots \\ v_t^{76} \end{bmatrix} = \boldsymbol{M^{dn}} \cdot \boldsymbol{V^{dn}} \tag{3}$$

Note that $\boldsymbol{M^{up}}$ and $\boldsymbol{M^{dn}}$ can be constructed by $\boldsymbol{M}$ if we view constraints row-wise and account for the relationship between configurations and their feasible transitions. That is, only $\boldsymbol{M}$ is needed to build constraints (1)-(3). In addition, because at any time interval $t$ only one configuration status $u_t^y$ is chosen, feasible transition $v_t^{yy'}$ can be uniquely determined accordingly. Thus, after mapping relationships $\boldsymbol{M}$, $\boldsymbol{M^{up}}$, and $\boldsymbol{M^{dn}}$ are constructed, PT-related variables $i_t^x$, $su_t^x$, and $sd_t^x$ could be fully represented by configuration-related variables as in (4), where sets $ON^x$, $UT^x$, and $DT^x$ can be respectively derived from corresponding rows of $\boldsymbol{M}$, $\boldsymbol{M^{up}}$, and $\boldsymbol{M^{dn}}$.

$$i_t^x = \sum_{y \in ON^x} u_t^y \,;\, su_t^x = \sum_{yy' \in UT^x} v_t^{yy'} \,;\, sd_t^x = \sum_{yy' \in DT^x} v_t^{yy'} \tag{4}$$

2) Aggregate Configuration Model

In aggregate configuration model, configurations with the same numbers of CTs and STs are aggregated as one to reduce computational complexity and provide MPs more operational flexibility. Using the same "2×1" example, 5 configurations and 12 transitions (6 upward and 6 downward) are needed in the aggregate configuration model, by combining configurations $CT_1$ and $CT_2$ as well as $CT_1$+ST and $CT_2$+ST in the complete configuration model. In turn, valid configurations are all off, 1CT, 2CT, 1CT+1ST and 2CT+1ST, denoted as $u_t^1$ to $u_t^5$. By further considering $CT_1$ having a higher priority to be turned on than $CT_2$, constraint (5) is derived for describing the mapping relationship between PTs and configurations. Similarly, constrains (6)-(7) describe the mappings between startup/shutdown and feasible transitions.

$$\begin{bmatrix} i_t^{c1} \\ i_t^{c2} \\ i_t^{s1} \end{bmatrix} = \begin{bmatrix} 0 & 1 & 1 & 1 & 1 \\ 0 & 0 & 1 & 0 & 1 \\ 0 & 0 & 0 & 1 & 1 \end{bmatrix} \begin{bmatrix} u_t^1 \\ \vdots \\ u_t^5 \end{bmatrix} \tag{5}$$

Although priorities of alternative CTs/STs are enforced in the aggregate configuration model, physical plant operators do have the flexibility to switch from one CT to another as long as it does not violate other physical limits (especially the features related to warmth status and minimum up/down time constraints). This would be valuable when one of alternative PTs encounters physical failures.

$$\begin{bmatrix} su_t^{c1} \\ su_t^{c2} \\ su_t^{s1} \end{bmatrix} = \begin{bmatrix} 1 & 1 & 0 & 0 & 0 & 0 \\ 0 & 1 & 1 & 0 & 0 & 1 \\ 0 & 0 & 0 & 1 & 1 & 0 \end{bmatrix} \begin{bmatrix} v_t^{12} \\ \vdots \\ v_t^{45} \end{bmatrix} \quad (6)$$

$$\begin{bmatrix} sd_t^{c1} \\ sd_t^{c2} \\ sd_t^{s1} \end{bmatrix} = \begin{bmatrix} 1 & 1 & 0 & 0 & 0 & 0 \\ 0 & 1 & 1 & 0 & 0 & 1 \\ 0 & 0 & 0 & 1 & 1 & 0 \end{bmatrix} \begin{bmatrix} v_t^{21} \\ \vdots \\ v_t^{54} \end{bmatrix} \quad (7)$$

3) Remarks

As discussed above and shown in Fig. 2, there are two sets of elements in CCGT models, i.e., PTs and configurations. CPBM focuses on the left subgraph of Fig. 2 with PT-related commitment, startup, and shutdown indicators, while transition logics are hard to describe with PT-related variables only. In comparison, CFBM focuses on the right subgraph of Fig. 2 with configuration-related commitment and transition indicators, in which PT properties are largely neglected.

In fact, these two sets of elements can be distinctly mapped with each other. Specifically, as valid configurations and feasible transitions are both unique at a time interval, PTs' status can be fully tracked with the established mapping relationship. That is, variables in CFBM contain sufficient information to derive commitment status of individual PTs. Reversely, $u_t^y$ and $v_t^{yy'}$ may also be linearly expressed via PT variables. However, because PT-related variables are not exclusive at any time interval like configuration-related variables, the binary nature of configuration-related variables may not hold via a linear mapping. That is, $u_t^y$ and $v_t^{yy'}$ could be integers instead of binaries, which deviate from their physical meanings. Actually, the range of integer values corresponding to a PT (when the PT is on) defines the feasible transition information in CPBM (i.e., constraints (20)-(30) and required input information in [5]).

In summary, with established mapping relationships, variables in CFBM can be used to indicate statuses of PTs and no additional variables are needed to depict individual PTs' features. Therefore, the proposed hybrid model remains the tight and compact fashion as discussed in [23], [27]-[28].

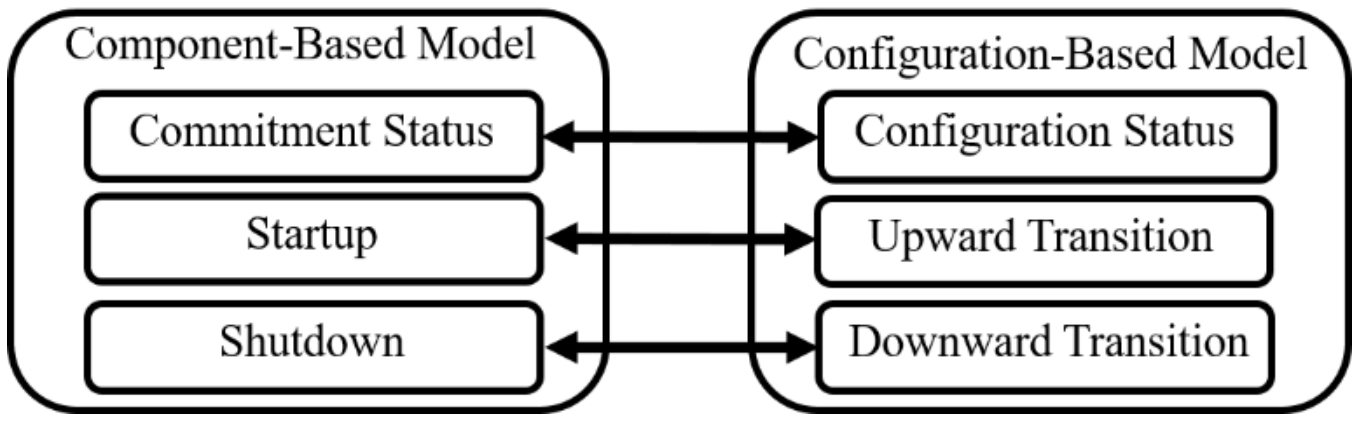

Fig. 2 Relationship between CPBM and CFBM

*B. Constraints with Configuration Features*

1) Incremental Energy Cost and No-load Cost

Operation cost of each CCGT, including energy cost $C^{P,y}(p_t^y)$ and no-load cost $C^{NL,y}$, is modeled via configuration status as shown in (8). Configuration based offers could be more accurate to reflect operation features of CTs/STs. That is, as output of an ST is strongly related to that of a CT, energy cost of the entire CCGT cannot be simply regarded as a direct summation of CTs and STs while ignoring their dependency. Here, locally ideal formulation with special ordered set of type 2 (SOS2) is used to piecewise linearize the energy cost function $C^{P,y}(p_t^y)$ [29], which has been proved to be the convex envelope of piecewise linear cost function and test results on the MISO system have shown great computational benefits [10].

$$\sum_{t\in T}\sum_{y\in Y}\left[C^{P,y}(p_t^y) + C^{NL,y}u_t^y\right] \quad (8)$$

2) Exclusion of Configurations

For each time interval $t$, all valid configurations of a CCGT are exclusive, as expressed in (9).

$$\sum_{y\in Y} u_t^y = 1 \qquad \forall t \quad (9)$$

3) Transition between Configurations

Similar to startup and shutdown variables of regular units, $\sum_{y'\in F^{*,y}} v_t^{y'y} / \sum_{y'\in F^{y,*}} v_t^{yy'}$, which refers to all feasible transitions to/from configuration $y$ at $t$, can be regarded as startup/ shutdown transitions of configuration $y$ at $t$. Thus, the logical relationship of commitment status for configuration $y$ between two consecutive time intervals can be expressed as in (10).

$$u_t^y - u_{t-1}^y = \sum_{y'\in F^{*,y}} v_t^{y'y} - \sum_{y'\in F^{y,*}} v_t^{yy'} \qquad \forall t, y\in Y \quad (10)$$

4) Capacity Limit

Each configuration of a CCGT, similar to regular units, has minimum/maximum capacity limit $\underline{P}_t^y / \overline{P}_t^y$.

$$u_t^y \cdot \underline{P}_t^y \le p_t^y \le u_t^y \cdot \overline{P}_t^y \qquad \forall t, y\in Y \quad (11)$$

5) Ramping Constraints

As different combinations of CTs and STs may have distinct ramp up/down rates, ramp rate is considered as a feature more closely associated with configuration status as shown in (12)-(13). In (12), if there is no transition to configuration $y$, i.e., $\sum_{y'\in F^{*,y}} v_t^{y'y} = 0$, ramp up constraint is ruled by $p_t^y - p_{t-1}^y \le RU^y$; if a CCGT transits to configuration $y$, i.e., $\sum_{y'\in F^{*,y}} v_t^{y'y} = 1$, (12) becomes $p_t^y - p_{t-1}^y \le \sum_{y'\in F^{*,y}} v_t^{y'y}\, RU0^y$. Same implications apply to ramp down constraint (13). Indeed, ramp constraints (12)-(13) has better performance in the MISO system with asymmetric startup/shutdown features, to be more specific, shutdown cost is not counted and shutdown ramping is not required in the MISO system. For instance, (13) implies that a resource can move from any $p_{t-1}^y \le \overline{P}_{t-1}^y$ to $p_t^y = 0$ in a shutdown period $t$.

$$p_t^y - p_{t-1}^y \le RU^y + \sum_{y'\in F^{*,y}} v_t^{y'y}\,(RU0^y - RU^y) \qquad \forall t, y\in Y \quad (12)$$

$$p_{t-1}^y - p_t^y \le u_t^y RD^y + (1-u_t^y)\overline{P}_{t-1}^y \qquad \forall t, y\in Y \quad (13)$$

*C. Constraints with Individual PT Features*

1) Minimum Up/Down Time Limits

Minimum up/down time of individual PTs are expressed as in (14). Same as regular units, [30] showed that constraints (14)-(15) is a facet projection of min-up/down constraints in the space of $(i^x, su^x, sd^x)$. Replacing PT related variables by (4), (14)-(15) can be equivalently expressed as in (16)-(17).

$$\sum_{t'=t-TU^x+1}^{t} su_{t'}^x \le i_t^x \qquad \forall t\in[TU^x, T], x\in X \quad (14)$$

$$\sum_{t'=t-TD^x+1}^{t} sd_{t'}^x \le 1 - i_t^x \qquad \forall t\in[TD^x, T], x\in X \quad (15)$$

$$\sum_{t'=t-TU^x+1}^{t}\sum_{yy'\in UT^x} v_{t'}^{yy'} \le \sum_{y\in ON^x} u_t^y \qquad \forall t\in[TU^x, T], x\in X \quad (16)$$

$$\sum_{t'=t-TD^x+1}^{t}\sum_{yy'\in DT^x} v_{t'}^{yy'} \le 1 - \sum_{y\in ON^x} u_t^y$$

$$\forall t \in [TD^x, T], x \in X \quad (17)$$

As constraints (16)-(17) of individual PTs are described independently, they could be added separately to alleviate computational burden if a certain PT is dominant or multiple CTs/STs within a same configuration present same features in terms of minimum up/down times. This will be applied and tested in the MISO system study in Section IV.C.

In addition, similar to regular units, the required initial online time $TU^{R,x}$ and offline time $TD^{R,x}$ of each PT due to the minimum up/down time requirement are expressed as in (18)-(19). Corresponding configuration commitment variables have to be fixed in those initial periods (20).

$$TU^{R,x} = \max\{0, (TU^x - TU^{0,x}) \cdot i^{0,x}\} \qquad \forall x \in X \quad (18)$$

$$TD^{R,x} = \max\{0, (TD^x - TD^{0,x}) \cdot (1 - i^{0,x})\} \qquad \forall x \in X \quad (19)$$

$$\sum_{y \in ON^x} u_t^y = i^{0,x} \qquad \forall t \in [1, TU^{R,x} + TD^{R,x}], x \in X \quad (20)$$

2) Time-Dependent Startup Cost

Time-dependent startup cost is usually applied on regular units [26] [31], which is related to unit's down time. Up to three startup types are modeled in the MISO system, including hot, intermediate, and cold. After discussing with MPs and plant operators, this feature of CCGTs is also strongly associated with individual PTs and should be more accurately described via PTs' startup cost instead of transition cost between configurations as in CFBM. Constraints (21)-(22) are used to determine startup types. Note that $T_1^{x,SU} = TD^x$ represents the hot startup type indicator. Similarly, $su_t^x$ and $sd_t^x$ could be substituted via constraint (4), and (21)-(22) can be equivalently expressed as in (23)-(24). In turn, startup cost of a CCGT is described as $\sum_{x \in X} \sum_{t \in T} \sum_{w \in W} C_w^{SU,x} \cdot \delta_{t,w}^x$.

$$\delta_{t,w}^x \le \sum_{t'=T_w^{x,SU}}^{T_{w+1}^{x,SU}-1} sd_{t-t'}^x, \quad \forall t \in \left[T_{w+1}^{x,SU}, T\right], w = \{1,2\}, x \in X \quad (21)$$

$$\sum_{w \in W} \delta_{t,w}^x = su_t^x \qquad \forall\, t, x \in X \quad (22)$$

$$\delta_{t,w}^x \le \sum_{t'=T_w^{x,SU}}^{T_{w+1}^{x,SU}-1} \sum_{yy' \in DT^x} v_{t-t'}^{yy'}$$
$$\forall t \in \left[T_{w+1}^{x,SU}, T\right], w = \{1,2\}, x \in X \quad (23)$$

$$\sum_{w \in W} \delta_{t,w}^x = \sum_{yy' \in UT^x} v_t^{yy'} \qquad \forall\, t, x \in X \quad (24)$$

In addition, for a PT $x$ having been initially off for $TD^{0,x}$ time slots, the initial startup type for the period of initial time slots within the scheduling horizon can be fixed as in (25).

$$\delta_{t,w}^x = 0 \qquad \forall t \in \left[T_{w+1}^{x,SU} - TD^{0,x}, T_{w+1}^{x,SU}\right], w = \{1,2\}, x \in X \quad (25)$$

Furthermore, for transitions involving startup of multiple PTs, e.g., from "all off" to 2 CTs, startup costs of individual PTs can be calculated separately since their on/off time durations could be different. In fact, for the "2×1" example in Fig. 1, up to 9 types of transition costs could be modeled with respect to 3 types of CT startup costs, which provides more accurate modelling capability than the 3 transition costs considered in CFBM [10].

3) Additional Specific Constraints

With the capability of modeling properties of individual PTs, certain specific constraints in the MISO system, such as maximum daily start time and energy limits that are enforced on configurations in CFBM [10], can be applied on individual PTs to better match the practical operation concerns.

### *D. Computational Complexity Comparison with CFBM*

This section compares the size of formulated problems and potential computational impacts of the proposed hybrid model and CFBM. Constraints (16)-(20) and (23)-(25) represent major difference in mathematical formulation between the hybrid model and CFBM presented in [10] and [23].

- *Number of Binary Variables*: with the mapping relationship (4) between configurations and PTs erected, no additional binary variables are required for representing min up/down behavior of individual PTs. As for time-dependent startup indicators, because the number of PTs is usually less than that of upward transitions, the hybrid model would encounter fewer binary variables and in turn less computational effort. For instance, the aggregated "2×1" example includes 3 PTs v.s. 6 upward transitions, which results in 15 binary variables (9 startup indicators and 6 upward transition variables) in the hybrid model v.s. 18 transition variables (6 upward transitions×3 types of startup) in CFBM.
- *Number of Constraints*: numbers of constraints for minimum up/down time and time-dependent startup cost are proportional to the number of PTs in the hybrid model, whereas proportional to the number of configurations in CFBM. As the number of PTs is usually smaller than that of configurations, e.g., 3 PTs v.s. 5 configurations in aggregated "2×1" example, this indicates a problem with fewer constraints.
- *Number of Non-zeros:* as behavior of individual PTs is indicated by the summation of configuration-based variables, more non-zeros are introduced, which could impact the tightness of formulation [32]-[33]. Compared with CFBM, (16) introduces more negative coefficients on LHS, which indicates a less tight formulation, while in (17) more positive coefficients are introduced on LHS, which indicates a tighter formulation. As hybrid model and CFBM deploy different offers, i.e., from PTs and configurations respectively, it might be difficult to accurately and fairly measure their integrated impacts on formulation tightness.

## III. Choice of Tight Formulations

In this section, we discuss the choice of formulations in the hybrid model for achieving better computational performance. Inspired by state-of-the-art research on formulation techniques of general UC and CCGT modeling [23], [27]-[31], [34]-[36], as well as our extensive computational studies on the MISO system, the following issues of the proposed hybrid CCGT model are discussed.

### *A. Feasible Transition*

To formulate feasible transitions of configurations, constraint (10) is enforced with logics on configurations, while it could also be enforced with PT logics as in (26) and further expressed with configuration-related variables as in (27).

$$i_t^x - i_{t-1}^x = su_t^x - sd_t^x \qquad \forall t, x \in X \quad (26)$$

$$\sum_{y \in ON^x} u_t^y - \sum_{y \in ON^x} u_{t-1}^y = \sum_{yy' \in UT^x} v_t^{yy'} - \sum_{yy' \in DT^x} v_t^{yy'}$$
$$\forall t, x \in X \quad (27)$$

Compared with (10), constraint (27) involves fewer constraints as the number of PTs is less than that of configurations, while more non-zeros $\sum_{y \in ON^x} u_t^y$ in each constraint instead of a single $u_t^y$ could potentially weaken the tightness of equality constraint. Computational performance of (27) and (10) will be compared in Section IV.B.3).

### *B. Choice of Binary/Continuous Variables*

$v_t^{yy'}$ represents the transition between configurations. As shown in (10), similar to startup/shutdown variables in UC, $v_t^{yy'}$ could be enforced as binary even defined as continuous variables. However, test results between 3-bin and 1-bin variable formulations on UC [34]-[35] shows that 3-bin could have better computational performance, because additional binary variables would induce better branching choices or the improved heuristic approach in presolving process of solvers. Computational performance of declaring $v_t^{yy'}$ as binary and continuous variables will be compared in Section IV.A.

The binary/continuous choice of time-dependent startup indicators would yield same property as suggested in [27] [31], as they are limited by shutdown indicator and impacted by offline times. Related computational performance will be compared in Section IV.B.

### *C. Exclusion of Transitions*

For each CCGT, at most one transition is allowed at each $t$ (28). Although (28) can be derived from (9)-(10), in most cases, including (28) in the formulation could significantly shorten the solving time by preventing unnecessary search branches and providing a tighter lower bound in hybrid model.

$$\sum_{y \in Y} \sum_{y' \in F^{y,*}} v_t^{yy'} \leq 1 \qquad \forall\, t \quad (28)$$

### *D. Reducing the Number of Time-Dependent Startup Types*

The number of time-dependent startup indicators could be reduced if startup costs of different startup types are the same. Specifically, for PTs with equal intermediate and hot startup costs $\delta_{t,1}^x$ can be removed, while $\delta_{t,2}^x$ can be removed for transitions with equal cold and intermediate startup costs. This technique is also applicable to regular units for reducing the number of startup type indicators.

### *E. Duct Burners (DB) Modelling*

DB is a common technology used by CCGT owners to enhance energy utilization and increase outputs. DB uses supplementary firing to increase heat energy of a gas turbine's exhaust, in an attempt to increase the output of a downstream HRSG and finally result in an extended ST capacity. Data from MPs in MISO shows that 14 out of 31 (45.2%) CCGTs are equipped with DBs. In addition, we also notice some unusual spikes in the last segments of certain configuration energy offers, which indicate that some MPs may have DBs combined with regular configurations, and the number of CCGTs actually equipped with DBs might be even higher.

Indeed, DB can be modeled as an extra PT in the hybrid model. Table II displays a typical transition matrix of a CCGT with DB. Other power enhancement approaches discussed in [3] [5] could be handled similarly, i.e., by adding additional configurations and virtual PTs.

TABLE II TYPICAL TRANSITION MATRIX OF A CCGT WITH DB

| | 1CT1ST | 1CT1ST1DB | 2CT1ST | 2CT1ST1DB | All Off |
|---|---|---|---|---|---|
| 1CT1ST | 1 | 1 | 1 | 0 | 1 |
| 1CT1ST1DB | 1 | 1 | 0 | 0 | 0 |
| 2CT1ST | 1 | 0 | 1 | 1 | 1 |
| 2CT1ST1DB | 0 | 0 | 1 | 1 | 0 |
| All Off | 1 | 0 | 1 | 0 | 1 |

## IV. CASE STUDIES

### *A. The Modified IEEE 118-Bus System*

The modified IEEE 118-bus system is studied to compare performance of the proposed hybrid model and CFBM in [23] with respect to different min up/down time limits. CFBM in [23] has been proved to be a tight and computational efficient CFBM formulation. 8 CCGTs with the aggregate "$2 \times 1$" configuration are studied, and their initial statuses are all set to offline for sufficient long period. Cost functions of regular units and CCGTs are piecewise linearized into 5 segments. The following four cases are studied, which are solved by Gurobi 6.5.2 on a personal laptop with Intel Core i7 2.6 GHz CPU and 12 GB memory.

- *Cases 1.1/1.2*: CFBM with configuration transition costs and min-up and down times of 2/3 hours for all configurations.
- *Case 2.1/2.2*, Hybrid model with turbine startup/shutdown costs and min-up and down times of 2/3 hours for all PTs.

#### 1) Economic Comparison with CFBM

Commitment status of the CCGT at Bus 49 in all four cases is illustrated in Table III. Specifically, in order to transit from offline to the most economic Configuration 5, in Cases 1.1 and 1.2, the CCGT needs to stay in Configurations 2 and 4 for 2 and 3 hours respectively, while all transitions only involve startup of PTs which would not incur any min-up/down violations from turbine's perspective. In comparison, both scenarios in Case 2 have the same commitment result with respect to different PT min-up/down time limits, which shows the effectiveness of the proposed hybrid model that can take advantage of accurate formulation on physical characters of individual PTs to identify the most cost efficient schedule. Comparison on costs of the two models in Table IV shows production cost saving under the hybrid model. In addition, for Case 2, as only startup is involved and no PT is shutdown till the end of horizon for all CCGTs, both scenarios in Case 2 have the same commitment result and therefore the same cost.

Indeed, with the capability of modeling min-up/down time limits for individual PTs, the proposed hybrid model could not only achieve lower cost, but also provide flexibility for both the system and MPs. From the system's point of view, as a CCGT may only involve consecutive upward/downward transitions due to its fast response characteristic for adapting hourly load fluctuations, the hybrid model could avoid setting the CCGT in a certain configuration for unnecessary long time to meet the configuration up/down time limits. It will also allow MPs to offer the true costs and operational parameters.

TABLE III CONFIGURATION STATUS OF THE CCGT AT BUS 49

| Mode | Hours (1-24) |
|---|---|
| Case 1.1 | 1 1 2 2 4 4 4 5 5 5 5 5 5 5 5 5 5 5 5 5 5 5 5 5 |
| Case 1.2 | **2 2** 2 **4** 4 4 5 5 5 5 5 5 5 5 5 5 5 5 5 5 5 5 5 5 |
| Case 2.1/2.2 | 1 1 **1 1 2** 4 **4** 5 5 5 5 5 5 5 5 5 5 5 5 5 5 5 5 5 |

TABLE IV COST COMPARISON ON THE MODIFIED 118-BUS SYSTEM

| | Case1.1 | Case1.2 | Case2.1 | Case2.2 |
|---|---|---|---|---|
| Obj ($) | 1,511,836.81 | 1,514,848.79 | 1,508,854.85 | 1,508,854.85 |

#### 2) Computational Performance Comparison with CFBM

Computational performance of Cases 1 and 2 is listed in Table V with three different optimality gaps, i.e., 0.5%, 0.1%, and 0.01%. Computational results from Cplex 12.6.3 are also reported for comparison. As stated in Section II.D, since different formats of offers adopted in the two models may impact their computational performance, while it is hard to equivalently convert information from one to the other, a general conclusion on computational performance of the two models may not be easily derived. Nevertheless, it is observed in both models that different min-up/down times could greatly impact the selection on binary variables and the presolving procedure in a solver [30] [37], which consequently present heterogeneous computational performance. Through extensive experiments, we observe that different load profiles, cost functions, choice of solvers, and optimality gap thresholds could all significantly impact the computational performance.

TABLE V COMPUTATIONAL TIME OF CFBM AND THE HYBRID MODEL (S)

| Solver | Gurobi | | | Cplex | | |
|---|---|---|---|---|---|---|
| Gap | 0.5% | 0.1% | 0.01% | 0.5% | 0.1% | 0.01% |
| Case1.1 | 57 | 68 | 151 | 9 | 25 | 173 |
| Case1.2 | 50 | 65 | 73 | 11 | 33 | 38 |
| Case2.1 | 57 | 65 | 97 | 11 | 32 | 163 |
| Case2.2 | 69 | 75 | 105 | 8 | 32 | 233 |

3) Comparing Different Formulations in the Hybrid Model

Different formulations of the hybrid model discussed in Section III A-C are compared. Specifically, F-1 denotes the formulation with constraint (28), F-2 denotes F-1 with constraint (27) instead of (10), F-3 refers to F-1 by adding constraint (27), F-4 represents F-1 with continuous transition variables, and F-5 is the formulation without constraint (28).

Computational performances of Cases 2.1-2.2 with different formulations are summarized in Table VI, in which "R" indicates the time for root node relaxation, and "Nodes" denotes the number of nodes explored to reach the predefined optimality gap. As shown in Table VI, comparison between F-1 and F-2 shows that although providing same logics to enforce feasible transitions, equation (27), as a more compact (fewer constraints) and less tight (more non-zeros) formulation, could greatly impact computational performance. This could also be reflected in Table VII that F-2 derives a lower LP relaxation value and a larger integrality gap. Indeed, it takes less time for F-2 to solve the LP as indicated by the computing time of root node relaxation, which could result in a less total time to reach a relatively larger gap (0.5% and 0.1%), but for smaller gap the effort is significantly increased as indicated by the time and the number of explored nodes. F-3 further shows the efforts by adding both (27) and (10), which could improve the performance with less time and fewer explored nodes in most cases as compared to F-1. Setting transition variables as continuous in F-4 results a formulation with fewer binary variables, which could also save time for computing LPs and achieve a higher optimality gap even with more nodes explored, but could introduce cumbersome effort in searching for a better solution with a smaller gap. F-5 shows the similar trend by removing constraint (28) when reaching a smaller gap.

TABLE VI COMPARISON OF DIFFERENT FORMULATIONS

| | | R(s) | 0.5% | | 0.1% | | 0.01% | |
|---|---|---|---|---|---|---|---|---|
| | | | Time (s) | Nodes | Time (s) | Nodes | Time (s) | Nodes |
| Case 2.1 | F-1 | 27 | 57 | 27 | 65 | 286 | 97 | 2979 |
| | F-2 | 14 | 45 | 27 | 59 | 638 | 710 | 58484 |
| | F-3 | 28 | 65 | 15 | 72 | 255 | 91 | 2076 |
| | F-4 | 29 | 51 | 28 | 65 | 569 | 1090 | 51211 |
| | F-5 | 31 | 60 | 27 | 102 | 468 | 2195 | 154768 |
| Case 2.2 | F-1 | 29 | 69 | 28 | 75 | 158 | 105 | 2911 |
| | F-2 | 22 | 53 | 29 | 87 | 739 | 1396 | 93230 |
| | F-3 | 27 | 63 | 27 | 70 | 248 | 85 | 1440 |
| | F-4 | 25 | 47 | 29 | 61 | 628 | 677 | 50932 |
| | F-5 | 25 | 55 | 27 | 159 | 1415 | 1430 | 124831 |

TABLE VII COMPARISON ON COMPUTATIONAL PERFORMANCE OF F-1 AND F-2

| | F-1 | | | F-2 | | |
|---|---|---|---|---|---|---|
| | MIPobj | LPobj | Igap | MIPobj | LPobj | Igap |
| Case 2.1 | 1508854.85 | 1500698.77 | 0.5405 | 1508862.27 | 1497334.92 | 0.7640 |
| Case 2.2 | 1508854.85 | 1500698.77 | 0.5405 | 1508854.84 | 1497334.92 | 0.7635 |

In summary, a less tight formulation may save effort for deriving an integer feasible solution, but could bring risks of failing to clear the market within the desired optimality gap and time limit. This could all bring severe issues in practice [8]. Therefore, it is better to deploy constraint (10) instead of (27), and set transition variables as binary. In addition, constraint (28) is also suggested and formulated in the MISO system for achieving better computational performance.

## *B. The MISO Test Case*

To facilitate CFBM cost/benefit study [10], MISO has requested MPs to submit configuration-based offers in early 2017 [9]. Data of 31 CCGTs are collected, including 119 configurations consisted by 100 PTs (including DBs). In addition, some CCGTs require a CT accompanied with an ST due to the physical structure needs, that is, the ST is combined with the CT in the aggregated configuration which further results in 74 combined turbines. After that, additional PT data are requested for building the proposed hybrid model, including min up/down time and three types of time-dependent startup costs/times. Three daily load profiles S1-S3 are randomly selected for comparing costs and computational performance of the hybrid model and CFBM in [10]. Four cases of the hybrid model are studied, including HM-1 which is based on CFBM while enforcing min-up/down time constraints on individual PTs, HM-2 which represents the full hybrid formulation described in Section II, HM-3 which changes startup indicators of HM-2 to continuous, and HM-4 which adds (27) to HM-2. Modeling techniques discussed in Sections III.C and D are applied to all the four hybrid models and CFBM. All test cases are performed on Intel Xeon X5680 @3.33GHz RAM 32 GB on with AIMMS 4.33.2 and CPLEX 12.7. The required relative optimality MILP gap is 0.1%.

TABLE VIII TEST RESULTS ON THE MISO SYSTEM

| | Model | Obj ($) | Saving ($) | Gap (%) | Time (s) | R (s) |
|---|---|---|---|---|---|---|
| S1 | CFBM | 52441054 | N/A | 0.035 | 465 | 191 |
| | HM-1 | 52421267 | 19787 | 0.075 | 485 | 176 |
| | **HM-2** | **52409022** | **32032** | **0.080** | **478** | **187** |
| | HM-3 | 52415287 | 25767 | 0.098 | 536 | 169 |
| | HM-4 | 52406778 | 34276 | 0.062 | 465 | 177 |
| S2 | CFBM | 76581171 | N/A | 0.039 | 461 | 205 |
| | HM-1 | 76567416 | 13755 | 0.054 | 444 | 226 |
| | **HM-2** | **76585759** | **-4588** | **0.096** | **356** | **193** |
| | HM-3 | 76593555 | -12384 | 0.099 | 355 | 193 |
| | HM-4 | 76574754 | 6417 | 0.070 | 425 | 201 |
| S3 | CFBM | 49313006 | N/A | 0.07 | 615 | 208 |
| | HM-1 | 49270622 | 42384 | 0.047 | 700 | 223 |
| | **HM-2** | **49291085** | **21921** | **0.094** | **599** | **200** |
| | HM-3 | 49289208 | 23798 | 0.089 | 694 | 192 |
| | HM-4 | 49294244 | 18762 | 0.099 | 514 | 218 |

Test results of the MISO system shown in Table VIII indicate that, the hybrid model could provide significant cost savings if difference in MILP gap is also considered (the cost increase in S2 is due to a larger MILP gap difference when terminated). With the relative optimality MILP gap of 0.1%, all three load scenarios show that HM-1 consumes more time than the full model HM-2. HM-2 also shows reasonable computational benefits over CFBM, by observing that computational times of S2 and S3 are respectively reduced by 22.8% and 2.6%, while only a slightly time increase in S1. In addition, although continuous startup indicator formulation in HM-3 consumes less time in root node relaxation, it may suffer significant increase in total time (i.e., 12.1% in S1 and 15.9% in S3). Thus, it is better to define them as integers. HM-4 shows same trends as observed in the 118-bus system studies, and tradeoff between the root node relaxation time and total time consumption in different formulations is also revealed. Future work could include exploring additional formulation techniques and tightening constraints for the hybrid model to achieve better performance.

## V. Conclusions

This paper proposes a hybrid CCGT model, by enjoying market offers from both configurations and individual PTs, for reflecting physical operation features more accurately and operating CCGTs more economic efficiently. The proposed hybrid model is compared with CFBM, and various tightening formulations with PT/configuration logics are analyzed to further improve computational performance of the hybrid model. Numerical studies show promising benefits on achieving better social welfare. The next step is to work close with MPs for more practical needs and to further improve performance of the hybrid model with state-of-the-art solvers.